\documentclass[11pt,reqno]{amsart}

\usepackage{amssymb}
\usepackage{amsmath}
\usepackage{amsfonts} 
\usepackage{amsbsy}
\usepackage{amscd}
\usepackage{amsthm}
\usepackage{geometry}    
\usepackage{courier}
\usepackage{xcolor}

\usepackage{graphicx}
\usepackage{amssymb}
\usepackage[capitalize]{cleveref}
\usepackage[framemethod=tikz]{mdframed} 
\usepackage{circuitikz}
\usetikzlibrary{decorations.pathmorphing, patterns,shapes}

\usetikzlibrary{arrows}
\usetikzlibrary{patterns}
\usepackage{xcolor}
\usepackage{caption}
\usepackage{subcaption}

\newcommand{\Hidden}[1]{}

\newtheorem{theorem}{Theorem}[section]

\newtheorem{prop}[theorem]{Proposition}

\newtheorem{conjecture}[theorem]{Conjecture} 
 
\theoremstyle{definition}

\newtheorem{definition}[theorem]{Definition} 
 
\theoremstyle{remark}

\title{An Introductory Survey of Recursions in the Computation of Resistance Distance}
\author{Emily J. Evans}
\address{Brigham Young University}
\email{EJEvans@math.byu.edu}
\author{Russell Jay Hendel}
\address{Towson University}
\email{RHendel@Towson.Edu}

\begin{document}

\begin{abstract}  This paper presents an introduction and expository account of a beautiful, current, and active application of recursions to the computation of resistance distance.  Resistance distance, also referred to as effective resistance, is a well-known graph metric that arises naturally by considering a graph as an electrical circuit; heuristically resistance distance measures both the number of paths between two vertices in a graph and the cost of each path. This topic finds applications in a rich array of fields including  social, biological, ecological, and transportation networks, chemistry,
graph theory, numerical linear algebra, and engineering. A variety of methods are used in the field to determine resistance distance including recursive, mathematical, and graphical techniques. Sequences familiar to the readers of the Fibonacci Quarterly such as the Fibonacci and Lucas sequences appear quite often in results in the literature. Twenty five to forty years ago there were a handful of papers on resistance that appeared in the Fibonacci Quarterly and the Proceedings and recently papers on the subject have appeared again. It is hoped that this introductory expository account will interest readers of the Quarterly to renew interest in this current and active field.   
\end{abstract}

\maketitle

KEYWORDS:
\textit{
recursions, effective resistance, resistance distance,  families of matrices, determinants, toeplitz, linear--2 tree, linear 3--tree}

\section{Introduction}
This paper introduces resistance distance, also referred to as effective resistance, a well-known metric on graphs that considers a graph as a circuit with edges represented as resistors.  Resistance distance encapsulates the number of paths between two vertices in a graph and the length of each path. This topic finds applications in a rich array of fields including mathematical connectivity in social, biological, ecological, and transportation networks ~\cite{Applications}, chemistry~\cite{rdmatrix,carmona2014effective,Cinkir,KleinRandic, klein2002resistance,klein2004random, kem1,peng2017kirchhoff, yang2014comparison, yang2008kirchhoff}, graph theory~\cite{bapatdvi,BapatWheels, DEVRIENDT202224,littleswim, MarkK, Ghosh, klein1997graph, ZHOU20172864}, numerical linear algebra~\cite{SpielSparse}, and engineering~\cite{Barooah06grapheffective}.   

The mathematical techniques used to study resistance distance, heavily rely on recursions and their closed-form representations also known as their Binet forms. Between 25 and 40 years ago a handful of papers appeared in the Fibonacci Quarterly and the Proceedings studying various aspects of resistance distance and their underlying recursions \cite{Ferri, Lahr, Nodine, Risk}. Recently the authors have published two papers continuing the tradition of studying resistance distance in the Fibonacci Quarterly \cite{Sarajevo, InfiniteArray}. It is hoped that this expository account will increase interest in this growing field.

At first blush, it might appear that resistance distance would only interest those studying electric circuits. Indeed, the original technical definition of the resistance distance between two nodes in a given graph $G,$ calculates resistance distance by considering the graph as an electric circuit where each edge is represented by a resistor. Then, given any two nodes $i$ and $j,$ assume that one unit of current flows into node $i$ and one unit of current flows out of node $j$.  The potential difference $v_i - v_j$ between nodes $i$ and $j$ needed to maintain this current is the {\it resistance distance} between $i$ and $j$.  However, resistance distance is a metric on a graph~\cite{KleinRandic, klein2002resistance}. This viewpoint immediately widens the applicability of resistance distance from electrical circuit theory to any application where graph modeling and distance metrics are present. Hence, it applies to all networks, social, biological, and many more, and is very useful in measuring, for example, atomic energy levels in chemistry.

 


\section{Families of Graphs and Computational Methods}

Having introduced basic concepts and the diverse applications in Section 1, this section reviews the rich variety of families of graphs that have been studied. We also review several techniques used to determine resistance distance in graphs.

First, there are many particular families of graphs for which general resistance formulae have been obtained such as circulant~\cite{circulant}, corona~\cite{corona}, octogonal~\cite{octogonal}, regular~\cite{regular}, (almost) complete bipartite~\cite{bipartite}, Cayley \cite{cayley2, cayley}, some particular cubic graphs~\cite{cubic}, ring clique~\cite{ringclique}, straight linear 1 and 2 trees ~ \cite{bef}, Apollonian~\cite{apollonian}, flower \cite{flower}, 	Sierpinski~\cite{sierpinski}, prism~\cite{Cinkir2} and ladder~\cite{Cinkir} graphs and network families. As already indicated, besides the theoretical interest in such graphs, they have applications in many fields.

\begin{figure} 
 \begin{center}
\resizebox{4in}{!}{ \begin{circuitikz}
 
 \draw (0,0) to [R, *-*,l=$R_B$] (4,0);
  \draw (0,0) to [R, *-*,l=$R_C$] (2,3.46);
  \draw (2,3.46) to[R, *-*,l=$R_A$] (4,0)
  {[anchor=north east] (0,0) node {$N_1$} } {[anchor=north west] (4,0) node {$N_3$}} {[anchor=south]  (2,3.46) node {$N_2$}};
  
  \draw (6,0) to [R, *-*,l=$R_1$] (8,1.73);
  \draw (6,3.46) to [R, *-*,l=$R_2$] (8,1.73);
 \draw (10.5,1.73) to [R, *-*,l=$R_3$] (8,1.73)
  {[anchor=north east] (6,0) node {$N_1$} } {[anchor=north west] (10.5,1.73) node {$N_3$}} {[anchor=south]  (6,3.46) node {$N_2$}};
 
 \end{circuitikz}}
 \end{center}
 \caption{$\Delta$ and $Y$ circuits illustrating the 
 mesh-star and star-mesh transforms. $N$ and $R$ with subscripts label nodes and edge-resistances respectively.}\label{fig:dy}
 \end{figure}
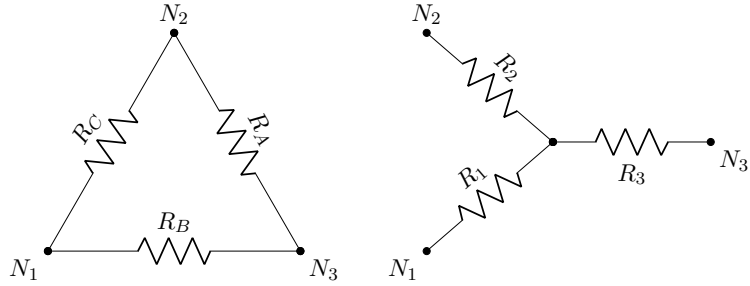

There are a variety of techniques used to compute resistance distance and we briefly review several here. A more complete summary of these techniques can be found in ~\cite{littleswim} which also includes many worked examples. Techniques can roughly be divided into two categories.  The first are circuit transformation techniques which use the properties of electrical circuits to transform a graph into a simpler form.  These techniques include the well-known series and parallel circuit transformations and the lesser-known mesh-star and star-mesh transformations~\cite[Section 2]{littleswim}.  Figure~\ref{fig:dy} illustrates the simplest star-mesh and mesh-star transform, known
as the Y-Delta and Delta-Y transformations. More formally these simple transformations are given by the following definitions.
\begin{definition}[$\Delta$--Y transformation]\label{def:dy}
Let $N_1, N_2, N_3$ be nodes and $R_A$, $R_B$ and $R_C$ be given resistances as shown in Figure~\ref{fig:dy}.  The transformed circuit in the ``Y'' format as shown in Figure~\ref{fig:dy} has the following resistances:
\begin{align*}
  R_1 &= \frac{R_BR_C}{R_A + R_B + R_C} \\
  R_2 &= \frac{R_AR_C}{R_A + R_B + R_C} \\
  R_3 &= \frac{R_AR_B}{R_A + R_B + R_C}.
\end{align*}
\end{definition}
\begin{definition}[Y--$\Delta$ transformation]\label{def:yd}
Let $N_1, N_2, N_3$ be nodes and $R_1$, $R_2$ and $R_3$ be given resistances as shown in Figure~\ref{fig:dy}.  The transformed circuit in the ``$\Delta$'' format as shown in Figure~\ref{fig:dy} has the following resistances:
\begin{align*}
  R_A &= \frac{R_1R_2 + R_2R_3 + R_1R_3}{R_1} \\
  R_B &= \frac{R_1R_2 + R_2R_3 + R_1R_3}{R_2} \\
  R_C &= \frac{R_1R_2 + R_2R_3 + R_1R_3}{R_3}.
\end{align*}
\end{definition} One downside of these techniques is that one must know beforehand the two vertices between which you want to determine the resistance distance.   Another downside is, that although it is always possible to apply circuit transformations to a planar graph, if the graph is not planar, one may not be able to transform a graph to a simpler form.

 The second category of methods is purely mathematical and consists of combinatorial, recursive, graph theoretic, numerical, and optimization techniques.  Several of the previously mentioned methods utilize the combinatorial Laplacian matrix.  We refer the interested reader to~\cite{littleswim} for worked examples of these techniques.  

When dealing with a family of graphs we utilize the notation 
$\{G_n\}_{n \ge 1},$ with the parameter $n$ typically indicating the number of vertices in the graph.  One question is whether a recursive relationship between maximal resistance distances in a family of graphs can be derived.  Since many of the techniques to determine resistance distance utilize the Laplacian matrix, and many techniques exist for recursive relationships in structured matrices (for example tridiagonal matrices~\cite{ELMIKKAWY2004669}, pentadiagonal matrices~\cite{pentadiag,recursion}, block tridiagonal matrices~\cite{molinari2008determinants}, and Toeplitz matrices~\cite{li2011calculating}), the idea of such a recursive relationship is not surprising.  

Towards this end, of deriving recursions, we formalize the definition of the Laplacian matrix.
\begin{definition}\label{def:Laplacian} Let $G$ be a simple, undirected graph on $n$ vertices enumerated in some order say $v_1, \dotsc, v_n.$ Define the $n \times n$ diagonal matrix $D,$ by $D_{i,i} = deg(v_i)$ with $D_{i,j}=0, \text{ if } i \neq j.$ Let $A,$ the adjacency matrix of G, be the $n \times n$ 0-1 matrix defined by
$A_{i,j} = 1, \text{ iff $v_i$ and $v_j$ are connected by an edge in $G$}.$ The (combinatorial) Laplacian, $L$ is an $n \times n$ matrix defined by $L = D-A.$ \end{definition}

 A general example, subsuming many others, is presented in \cite{recursion} which shows, that the determinants of the general pentadiagonal family of matrices governed by five parameters, satisfy a sixth-order recursion whose roots can be explicitly calculated.  These recursions allow computation of the closed Binet forms and, as a consequence, closed-formula for resistances. As \cite{recursion} shows, such an approach is computationally quicker than many other methods for computing resistance distance.

A fundamental result of ~\cite{bapatdvi} is that the resistance distance can be computed using determinants associated with the underlying  Laplacian with specific rows and columns deleted. The authors~\cite{Automated} have recently suggested an automated procedure based on this method mimicking the approach in~\cite{recursion}.

\section{Sample families of graphs and relationships to Fibonacci numbers}\label{sec:examples}
Having presented the ideas, applications, methods, and graph families associated with resistance distance,  
this section presents a handful of examples of families of graphs each example accompanied by a definition,  a figure, a representative Laplacian matrix, and  results. We omit proofs noting that several examples are subsumed under the penta-diagonal family of graphs which is thoroughly investigated in \cite{recursion}. We particularly note the presence of the Fibonacci and Lucas sequences in results suggesting a natural relationship between this field and the themes of the Fibonacci Quarterly.

\subsection{Path graphs}\label{subsec:path}   
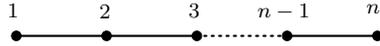
\begin{figure}
\begin{center}

\begin{tikzpicture}[line cap=round,line join=round,>=triangle 45,x=1.0cm,y=1.0cm, scale = 1.2]
\draw [line width=.8pt] (-3.,2.)-- (-2.,2.);
\draw [line width=.8pt] (-2.,2.)-- (-1.,2.);
\draw [line width=.8pt,dotted] (-1.,2.)-- (0.,2.);
\draw [line width=.8pt] (0.,2.)-- (1.,2.);
\begin{scriptsize}
\draw [fill=black] (-3.,2.) circle (1.5pt);
\draw[color=black] (-3.0357596415464108,2.283574198210263) node {$1$};
\draw [fill=black] (-2.,2.) circle (1.5pt);
\draw[color=black] (-2.0176831175804275,2.269627944457304) node {$2$};
\draw [fill=black] (-1.,2.) circle (1.5pt);
\draw[color=black] (-1.0274991011203625,2.283574198210263) node {$3$};
\draw [fill=black] (0.,2.) circle (1.5pt);
\draw[color=black] (-0.037315084660296774,2.269627944457304) node {$n-1$};
\draw [fill=black] (1.,2.) circle (1.5pt);
\draw[color=black] (0.9528689317997687,2.297520451963222) node {$n$};
\end{scriptsize}
\end{tikzpicture}

%
%
\end{center}
\caption{A path graph on $n$ vertices.}
\label{fig:path}
\end{figure}

The simplest family of graphs under consideration is the family of path graphs as shown in Figure~\ref{fig:path}.  The adjacency matrix for the path graph is the banded matrix with ones on both the super and sub diagonal, and hence the $n \times n$ Laplacian is given by
\[L_G=\left[\begin{array}{rrrrrrrr}
1 & -1 & 0 & 0  &   0& \dots & \dots &0\\
-1 & 2 & - 1 & 0 &    0&  0& \ddots& \vdots\\
0 &-1 & 2 & - 1 & 0 &    0&  \ddots&  \vdots\\
0 &0 &-1 & 2 & - 1 & 0 &  \ddots&  0\\
0 & \ddots &\ddots & \ddots &  \ddots &  \ddots&  \ddots & 0\\
\vdots & \ddots & 0 & 0 &  -1 &  2&  -1 & 0\\
\vdots & \ddots  &0 & 0 & 0 &  -1 &  2&  -1 \\[1.7mm]
0  &\dots& \dots & 0 &0& 0 &  -1 & 1
\end{array}\right].\]

The resistance distance between two vertices in the path graph is the same as the distance, that is, for a path graph as shown in Figure~\ref{fig:path}, with $n \ge 2,$ the resistance distance between node $1$ and node $n$ is $n-1$. 

The literature traditionally formulates resistance distances in terms of the familiar Fibonacci-Lucas sequences and linear functions of $n.$ However, it is equally possible to formulate results in terms of linear homogeneous recursions with constant coefficients. In this case, the sequence $G_n=n-1$ satisfies the recursion $G_{n}=2 G_{n-1} - G_{n-2},$ with   initial conditions $G_0=-1$ and $G_1=0.$ Like all recursive sequences, the sequence may be extended to all negative and positive integers, but its applicability to path graphs is for $n \ge 2.$

\subsection{Straight Linear 2--trees}\label{sub:2tree}
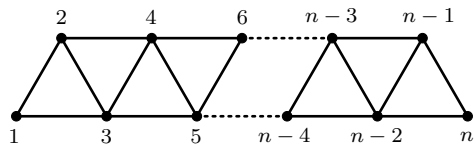
\begin{figure}[h!]
\begin{center}

\begin{tikzpicture}[line cap=round,line join=round,>=triangle 45,x=1.0cm,y=1.0cm,scale = 1.2]
\draw [line width=1.pt] (-3.,0.)-- (-2.,0.);
\draw [line width=1.pt] (-2.,0.)-- (-1.,0.);
\draw [line width=1.pt,dotted] (-1.,0.)-- (0.,0.);
\draw [line width=1.pt] (0.,0.)-- (1.,0.);
\draw [line width=1.pt] (1.,0.)-- (2.,0.);
\draw [line width=1.pt] (2.,0.)-- (1.5,0.866025403784435);
\draw [line width=1.pt] (1.5,0.866025403784435)-- (1.,0.);
\draw [line width=1.pt] (1.,0.)-- (0.5,0.8660254037844366);
\draw [line width=1.pt] (0.5,0.8660254037844366)-- (0.,0.);
\draw [line width=1.pt] (-0.5,0.8660254037844378)-- (-1.,0.);
\draw [line width=1.pt] (-1.,0.)-- (-1.5,0.8660254037844385);
\draw [line width=1.pt] (-1.5,0.8660254037844385)-- (-2.,0.);
\draw [line width=1.pt] (-2.,0.)-- (-2.5,0.8660254037844388);
\draw [line width=1.pt] (-2.5,0.8660254037844388)-- (-3.,0.);
\draw [line width=1.pt] (-2.5,0.8660254037844388)-- (-1.5,0.8660254037844385);
\draw [line width=1.pt] (-1.5,0.8660254037844385)-- (-0.5,0.8660254037844378);
\draw [line width=1.pt,dotted] (-0.5,0.8660254037844378)-- (0.5,0.8660254037844366);
\draw [line width=1.pt] (0.5,0.8660254037844366)-- (1.5,0.866025403784435);
\begin{scriptsize}
\draw [fill=black] (-3.,0.) circle (1.5pt);
\draw[color=black] (-3.02279181666165,-0.22431183338253265) node {$1$};
\draw [fill=black] (-2.,0.) circle (1.5pt);
\draw[color=black] (-2.0001954862580344,-0.22395896857501957) node {$3$};
\draw [fill=black] (-2.5,0.8660254037844388) circle (1.5pt);
\draw[color=black] (-2.5018465162673555,1.100526386601008717) node {$2$};
\draw [fill=black] (-1.5,0.8660254037844385) circle (1.5pt);
\draw[color=black] (-1.5081915914412003,1.100526386601008717) node {$4$};
\draw [fill=black] (-1.,0.) circle (1.5pt);
\draw[color=black] (-1.0065405614318794,-0.22290037415248035) node {$5$};
\draw [fill=black] (-0.5,0.8660254037844378) circle (1.5pt);
\draw[color=black] (-0.4952423962300715,1.100526386601008717) node {$6$};
\draw [fill=black] (0.,0.) circle (1.5pt);
\draw[color=black] (-0.03217990699069834,-0.22431183338253265) node {$n-4$};
\draw [fill=black] (0.5,0.8660254037844366) circle (1.5pt);
\draw[color=black] (0.4887653934035965,1.103344389715898) node {$n-3$};
\draw [fill=black] (1.,0.) circle (1.5pt);
\draw[color=black] (0.9904164234129174,-0.22431183338253265) node {$n-2$};
\draw [fill=black] (1.5,0.866025403784435) circle (1.5pt);
\draw[color=black] (1.5692445349621338,1.1042991524908385) node {$n-1$};
\draw [fill=black] (2.,0.) circle (1.5pt);
\draw[color=black] (1.993718483431559,-0.22431183338253265) node {$n$};
\end{scriptsize}
\end{tikzpicture}
\end{center}
\caption{A straight linear 2-tree}
\label{fig:2tree}
\end{figure}

The next family of graphs 
 under consideration, a  generalization of the path graph is the so-called straight linear two tree, sometimes referred to as a two-path in the literature.  This graph is shown in Figure \ref{fig:2tree};   the $n \times n$ Laplacian matrix is:
\[L_G=\left[\begin{array}{rrrrrrrr}
\label{equ:2tree}
2 & -1 & -1 & 0  &   0& \dots & \dots &0\\
-1 & 3 & - 1 & -1 &    0&  0& \ddots& \vdots\\
-1 &-1 & 4 & - 1 & -1 &    0&  \ddots&  \vdots\\
0 &-1 &-1 & 4 & - 1 & -1 &  \ddots&  0\\
0 & \ddots &\ddots & \ddots &  \ddots &  \ddots&  \ddots & 0\\
\vdots & \ddots & 0 & -1 &  -1 &  4&  -1 & -1\\
\vdots & \ddots  &0 & 0 & -1 &  -1 &  3&  -1 \\[1.7mm]
0  &\dots& \dots & 0 &0& -1 &  -1 & 2
\end{array}\right].\]

    The resistance distance between any two vertices of this graph is known (see~\cite{bef}) but we only give the result between the vertices of degree two.
	\begin{theorem} \cite{bef} \label{the:2tree}
 
	Let $G_n$ be the linear 2--tree with $n$ vertices with $n  \ge 2.$ Then the resistance distance between nodes $1$ and $n$ is given by
\begin{equation}\label{equ:theorem31}
	r(1,n)=\frac{2F_{n-1}^2}{L_{n-1}L_{n-2}}+\sum_{i=1}^{n-3}\frac{F_iF_{i+1}}{L_iL_{i+1}}
	 = \frac{n-1}{5} + \frac{4F_{n-1}}{5L_{n-1}},
\end{equation}
  where $F_k$ and $L_k$ refer to the $k$th Fibonacci and Lucas numbers respectively.
	\end{theorem}

In this case,  the numerator and denominator of the second summand of the  right-most member of \eqref{equ:theorem31}   both satisfy the (shifted) Fibonacci-Lucas recursion, $G_n = G_{n-1} + G_{n-2},$ while the first summand   satisfies the recursion $G_n = 2G_{n-1} - G_{n-2},$ with initial conditions $G_0=-1, G_1=0.$ Since sums of recursive sequences are themselves recursive it is theoretically possible to express $r(1,n)$
as a ratio of a numerator and denominator each satisfying a single recursion. But the order of those recursions would be greater than 2 and would not add additional insight. 

\subsection{Straight linear 3--trees}\label{sub:linear3tree}

\begin{figure}[ht!]
    \centering
\begin{tikzpicture}[line cap=round,line join=round,>=triangle 45,x=1.0cm,y=1.0cm, scale = .8]
\draw [line width=.8pt] (0.7119648044207935,0.056502145327751115)-- (-0.19048706476012342,2.0870188509848124);
\draw [line width=.8pt] (-0.19048706476012342,2.0870188509848124)-- (2.6636787078204596,1.3307124522628444);
\draw [line width=.8pt] (2.6636787078204596,1.3307124522628444)-- (0.7119648044207935,0.056502145327751115);
\draw [line width=.8pt] (0.7119648044207935,0.056502145327751115)-- (1.2007929002271236,2.481841543751463);
\draw [line width=.8pt] (1.2007929002271236,2.481841543751463)-- (-0.19048706476012342,2.0870188509848124);
\draw [line width=.8pt] (1.2007929002271236,2.481841543751463)-- (2.6636787078204596,1.3307124522628444);
\draw [line width=.8pt] (1.5392123511699674,3.685110702659351)-- (-0.19048706476012342,2.0870188509848124);
\draw [line width=.8pt] (1.5392123511699674,3.685110702659351)-- (1.2007929002271236,2.481841543751463);
\draw [line width=.8pt] (1.5392123511699674,3.685110702659351)-- (2.6636787078204596,1.3307124522628444);
\draw [line width=.8pt] (1.5392123511699674,3.685110702659351)-- (3.43812149257148,3.478298815972058);
\draw [line width=.8pt] (3.43812149257148,3.478298815972058)-- (2.6636787078204596,1.3307124522628444);
\draw [line width=.8pt] (1.2007929002271236,2.481841543751463)-- (3.43812149257148,3.478298815972058);
\draw [line width=.8pt] (1.2007929002271236,2.481841543751463)-- (4.2,2.48);
\draw [line width=.8pt] (1.5392123511699674,3.685110702659351)-- (4.2,2.48);
\draw [line width=.8pt] (3.43812149257148,3.478298815972058)-- (4.2,2.48);
\begin{scriptsize}
\draw [fill=black] (0.7119648044207935,0.056502145327751115) circle (2.pt);
\draw[color=black] (0.6743626432049219,-0.19731244287938157) node {1};
\draw [fill=black] (-0.19048706476012342,2.0870188509848124) circle (2.pt);
\draw[color=black] (-0.547707596310903,2.171623713720523) node {2};
\draw [fill=black] (2.6636787078204596,1.3307124522628444) circle (2.pt);
\draw[color=black] (2.836486913117535,1.2691718445396072) node {$3$};
\draw [fill=black] (1.2007929002271236,2.481841543751463) circle (2.pt);
\draw[color=black] (0.9563788523239585,2.773258293174467) node {$4$};
\draw [fill=black] (1.5392123511699674,3.685110702659351) circle (2.pt);
\draw[color=black] (1.3700026256985456,3.995328532690291) node {$5$};
\draw [fill=black] (3.43812149257148,3.478298815972058) circle (2.pt);
\draw[color=black] (3.6449333792587733,3.6757101623553834) node {$6$};
\draw [fill=black] (2.6636787078204605,1.3307124522628435) circle (2.0pt);

\draw [fill=black] (4.2007929002271236,2.481841543751463) circle (2.pt);
\draw[color=black] (4.4563788523239585,2.773258293174467) node {$7$};
\end{scriptsize}

\end{tikzpicture}

\caption{A straight linear 3--tree on $7$ vertices.  The vertex numbering is such that vertex 1 and vertex $7$ have degree 3.}
\label{fig:st3tree}
\end{figure}
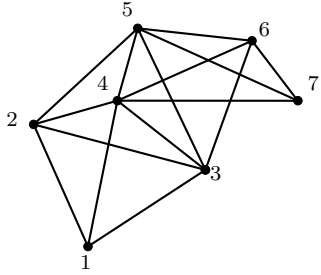

We can further generalize the path graph and 2-tree by considering the \textit{straight linear 3--tree} (sometimes referred to in the literature as a 3--path). The graph of the  3--tree has two vertices of degree three. Its  adjacency matrix consists of ones on the first three super and sub diagonals and zero elsewhere.  An example of such a tree on seven vertices is shown in Figure \ref{fig:st3tree}.

The $n \times n$ Laplacian matrix, $L$ of the straight linear 3--tree is given by

\begin{equation}\label{equ:PLinear3}
L_G=\begin{bmatrix}	
3 & -1 & -1 & -1 & \dotsc  &0 &0 &  0 & 0\\
-1 & 4 & -1 & -1 & \dotsc  &0 &0 &  0 & 0 \\ 
-1 & -1 & 5 & -1 &\dotsc   &0 & 0 & 0 & 0\\
-1 & -1 & -1 & 6 &\dotsc   &0 & 0 & 0 & 0  \\
\vdots & \ddots & \ddots & \ddots &\ddots &\ddots
\ddots & \ddots &  \vdots \\
0 & 0 & 0 & 0 & \dotsc &6  &-1  &-1 & -1  \\
0 & 0 & 0 & 0 & \dotsc &-1  &5  & -1 & -1  \\
0 & 0 & 0 & 0 & \dotsc &-1  &-1 &  4 & -1 \\
0 & 0 & 0 & 0 & \dotsc &-1  &-1  & -1 & 3\\
\end{bmatrix}.
\end{equation}
The middle of the matrix continues with sixes on the diagonal and negative ones on the first three super and subdiagonals.   The following limit has been conjectured.
\begin{conjecture}\label{con:onefourteenth} \cite{bef}  Let $G$ be the straight linear 3--tree,  with $n$ vertices and $H$ be the straight linear 3--tree with $n+1$ vertices. Let $r_G(1,n)$ and  ($r_H(1,n+1)$) indicate the total resistance between the two corner nodes 1 and $n$ ($1$ and $n+1)$. 
	Then 
	\[\lim_{n\rightarrow \infty} r_{H} (1, n+1) - r_G(1,n) = \frac{1}{14}.\]
 \end{conjecture}

Recently, the authors \cite{Automated} have, using the methods of recursions satisfied by a family of determinants, explicitly given the recursions of two sequences whose ratio is the resistance distance. The two sequences satisfy recursions of degrees 5 and 14, and have explicit roots lying in   quartic extensions of the rationals. 
    The denominator sequence satisfies the recursion
$G_n = 5 G_{n-1} - 3 G_{n-2} + 3 G_{n-3} -5 G_{n-4} + G_{n-5}$ with appropriate initial conditions while the numerator sequence satisfies the recursion
 $G_n = 7 G_{n-1} - 7 G_{n-2} - 98 G_{n-4} +56 G_{n-5} -56 G_{n-6} +198 G_{n-7} - 56 G_{n-8} + 56 G_{n-9} -98 G_{n-10} - 7 G_{n-12} + 7 G_{n-13}-G_{n-14}$ with appropriate initial conditions.

Their Binet forms have coefficients which can be explicitly presented and lie in quartic extensions of the rationals. Using these explicit forms coupled with the determinant formula for resistance distance, the conjecture can be proven. 
 
\subsection{Ladder graphs}\label{sub:ladder} 
\begin{figure}[ht!]
\begin{center}
\begin{tikzpicture}[line cap=round,line join=round,>=triangle 45,x=1.0cm,y=1.0cm,scale = 1]
\draw [line width=.8pt] (3.,4.)-- (3.,3.);
\draw [line width=.8pt] (3.,4.)-- (4.,4.);
\draw [line width=.8pt] (3.,3.)-- (4.,3.);
\draw [line width=.8pt] (4.,4.)-- (4.,3.);
\draw [line width=.8pt] (4.,4.)-- (5.,4.);
\draw [line width=.8pt] (4.,3.)-- (5.,3.);
\draw [line width=.8pt] (5.,4.)-- (5.,3.);
\draw [fill=black] (5.3,4.) circle (.6pt);
\draw [fill=black] (5.5,4.) circle (.6pt);
\draw [fill=black] (5.7,4.) circle (.6pt);
\draw [fill=black] (5.3,3.) circle (.6pt);
\draw [fill=black] (5.5,3.) circle (.6pt);
\draw [fill=black] (5.7,3.) circle (.6pt);
\draw [line width=.8pt] (6.,4.)-- (6.,3.);
\draw [line width=.8pt] (6.,4.)-- (7.,4.);
\draw [line width=.8pt] (6.,3.)-- (7.,3.);
\draw [line width=.8pt] (7.,4.)-- (7.,3.);
\begin{small}
\draw [fill=black] (3.,4.) circle (1.6pt);
\draw[color=black] (3,4.2) node {$1$};
\draw [fill=black] (3.,3.) circle (1.6pt);
\draw[color=black] (3,2.8) node {$2$};
\draw [fill=black] (4.,4.) circle (1.6pt);
\draw[color=black] (4,4.2) node {$3$};
\draw [fill=black] (4.,3.) circle (1.6pt);
\draw[color=black] (4,2.8) node {$4$};
\draw [fill=black] (5.,4.) circle (1.6pt);
\draw [fill=black] (5.,3.) circle (1.6pt);
\draw [fill=black] (6.,4.) circle (1.6pt);
\draw [fill=black] (6.,3.) circle (1.6pt);
\draw [fill=black] (7.,4.) circle (1.6pt);
\draw[color=black] (7,4.2) node {$2m-1$};
\draw [fill=black] (7.,3.) circle (1.6pt);
\end{small}
\begin{scriptsize}

\draw[color=black] (7,2.8) node {$2m$};

\end{scriptsize}
\end{tikzpicture}

\end{center}
\caption{The ladder graph on $n=2m$ vertices.}\label{fig:laddergraph}
\end{figure}
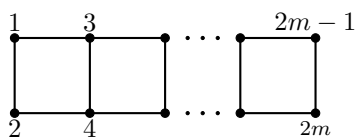

An alternative generalization of the path graph is the so-called ladder graphs on $n=2m$ vertices as illustrated in Figure~\ref{fig:laddergraph}. This   graph is the Cartesian product of $P_m$ and $P_2$. The first known resistance distance results were obtained by Cinkir~\cite{Cinkir}. 

In particular we have the following theorem:
\begin{theorem} \cite{Cinkir}
Let $G_n$ be the ladder graph with $n=2m$ vertices, labeled as in Figure~\ref{fig:laddergraph}. Define a generalized Fibonacci sequence by $H_0 =0,$   $H_1=1,$ $ H_n = 4 H_{n-1}- H_{n-2}, n \ge 2.$ Then the resistance distance between nodes $1$ and $2m,$ notationally indicated by $r(1,2m)$ is given by
	\[r(1,2m)=  -1 + \frac{H_{2m}}{2 H_m^2}. \] 
\end{theorem}

In this case, little would be gained by giving the recursions satisfied by the numerator and denominator (instead of the formulation in terms of squares and even indices of the recursion $H$).

The Laplacian matrix has an easily describable structure.  The first super and sub diagonal entries alternate between 0 and $-1$ and the second super and sub diagonals are identically $-1$ as shown in~\eqref{equ:ladderlaplacian}
\begin{equation}\label{equ:ladderlaplacian}L_G=\left[\begin{array}{rrrrrrrr}
2 & -1 & -1 & 0  &   0& \dots & \dots &0\\
-1 & 2 & 0 & -1 &    0&  0& \ddots& \vdots\\
-1 &0 & 3 & - 1 & -1 &    0&  \ddots&  \vdots\\
0 &-1 &-1 & 3 & 0 & -1 &  \ddots&  0\\
0 & \ddots &\ddots & \ddots &  \ddots &  \ddots&  \ddots & 0\\
\vdots & \ddots & 0 & -1 &  -1 &  3&  0 & -1\\
\vdots & \ddots  &0 & 0 & -1 &  0 &  2&  -1 \\[1.7mm]
0  &\dots& \dots & 0 &0& -1 &  -1 & 2
\end{array}\right].\end{equation}

 \subsection{Fan Graphs}\label{sub:fans}
 
\begin{figure}
\[\begin{array}{c}
\begin{tikzpicture}[line cap=round,line join=round,>=triangle 45,x=1.0cm,y=1.0cm,scale = 1.5]
\draw [line width=.8pt] (5.,3.)-- (4.191784354657935,3.848009081626107);
\draw [line width=.8pt] (4.39006772869402,1.9998427342084801)-- (3.8288722643919675,3.028138208643161);
\draw [line width=.8pt] (4.39006772869402,1.9998427342084801)-- (5.,3.);
\draw [line width=.8pt] (5.,3.)-- (3.8288722643919675,3.028138208643161);
\draw [line width=.8pt] (3.8288722643919675,3.028138208643161)-- (4.191784354657935,3.848009081626107);
\draw [line width=.8pt] (4.191784354657935,3.848009081626107)-- (5.,3.);
\draw [line width=.8pt] (5.,3.)-- (5.028138208643161,4.1711277356080325);
\draw [line width=.8pt] (5.028138208643161,4.1711277356080325)-- (4.191784354657935,3.848009081626107);
\draw [line width=.8pt] (5.,3.)-- (5.848009081626107,3.8082156453420652);
\draw [line width=.8pt] (5.848009081626107,3.8082156453420652)-- (6.171127735608033,2.9718617913568393);
\draw [line width=.8pt] (6.171127735608033,2.9718617913568393)-- (5.,3.);
\draw [line width=.8pt] (5.,3.)-- (5.5611954643020525,1.9717045255653174);
\draw [line width=.8pt] (5.5611954643020525,1.9717045255653174)-- (6.171127735608033,2.9718617913568393);
\begin{scriptsize}
\draw [fill=black] (5.,3.) circle (1.5pt);
\draw[color=black] (5,2.7) node {$1$};
\draw [fill=black] (4.191784354657935,3.848009081626107) circle (1.5pt);
\draw[color=black] (4.05,4.05) node {$4$};
\draw [fill=black] (3.8288722643919675,3.028138208643161) circle (1.5pt);
\draw[color=black] (3.5,3.08) node {$3$};
\draw [fill=black] (5.028138208643161,4.1711277356080325) circle (1.5pt);
\draw[color=black] (5,4.381775026148471) node {$5$};
\draw [fill=black] (5.848009081626107,3.8082156453420652) circle (1.5pt);
\draw[color=black] (6.1,4) node {$k-2$};
\draw [fill=black] (6.171127735608033,2.9718617913568393) circle (1.5pt);
\draw[color=black] (6.62,3.05) node {$k-1$};
\draw [fill=black] (4.39006772869402,1.9998427342084801) circle (1.5pt);
\draw [fill=black] (5.5611954643020525,1.9717045255653174) circle (1.5pt);
\draw[color=black] (5.25,1.8) node {$k$};
\draw[color=black] (4.2,1.8) node {$2$};
\draw [fill=black] (5.438073645134634,3.989671690475049) circle (.9pt);
\draw [fill=black] (5.233105926888897,4.08039971304154) circle (.9pt);
\draw [fill=black] (5.64304136338037,3.898943667908557) circle (.9pt);
\end{scriptsize}
\end{tikzpicture}
\end{array}\]
\caption{The fan graph on $k$ vertices}\label{fig:fan}
\end{figure}
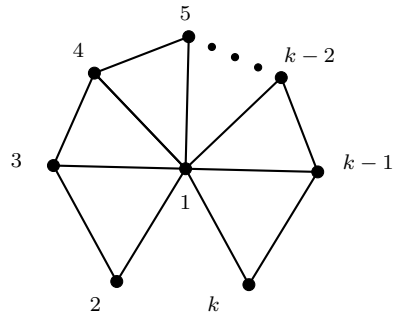

Another possible generalization of the path graph is obtained by joining the path graph with a singleton vertex.  This results in the so-called fan graph, which is shown in Figure~\ref{fig:fan}.

The Laplacian of the fan graph is given by
\[L_G=\left[\begin{array}{rrrrrrrr}
2 & -1 & 0 & 0  &  \ldots& \dots & 0 &-1\\
-1 & 3 & - 1 & 0 &    0&  0& \ddots& \vdots\\
0 &-1 & 3 & - 1 & 0 &    0&  \ddots&  \vdots\\
0 &0 &-1 & 3 & - 1 & 0 &  \ddots&  -1\\
\vdots & \ddots &\ddots & \ddots &  \ddots &  \ddots&  \ddots & -1\\
\vdots & \ddots & 0 & 0 &  -1 &  3&  -1 & -1\\
0 & \ddots  &0 & 0 & 0 &  -1 &  2&  -1 \\[1.7mm]
-1  &\dots& \dots & -1 &-1& -1 &  -1 & k-1
\end{array}\right].\]
 
 
Notice, that in contrast to the prior examples, the Laplacian is not banded.   Like many of the examples listed prior, a closed formula for the resistance distance between node $1$ and node $n$ is known.
\begin{prop}\cite{BapatWheels}\label{pro:bapatfans}
Let $k \geq 2$ be a positive integer. Then for $i= 1, \ldots, k-1$, the resistance distance between node $i$ and node $k$ in the fan graph is given by 
\begin{gather*}
r(i,k) = \frac{F_{2(k-1-i)+1}F_{2i-1}}{F_{2k-2}},
\end{gather*}
where $F_i$ is the $i$th Fibonacci number.
\end{prop}

Alternatively, we may formulate $r(i,k)$ in terms of recursive sequences   as follows. $r(i,k)=\frac{N(i,k)}{D_k},$ with for each fixed $i,$ $N(i,k)$ satisfies the recursion
$G_n=3 G_{n-1}-G_{n-2},$ with initial conditions 
$\{F_{2i-1} F_{1-2(i+1)}, 
F_{2i-1} F_{1-2i},
F_{2i-1} F_{1+2(2-i)} \}$, and with 
$(D_k)_{k \ge 2}$ satisfying the recursion  $G_n=3 G_{n-1}-G_{n-2}$ with  initial conditions, 
$\{G_0, G_1, G_2\} = \{-1,0,1\}.$

\subsection{Wheel Graphs}\label{sub:wheels}
An easy modification to the fan graph is to add an edge between node $1$ and node $n-1$.  This results in the so-called wheel graph as presented in Figure \ref{fig:wheels}.

The Laplacian of the wheel graph is given by
\[L_G=\left[\begin{array}{rrrrrrrr}
3 & -1 & 0 & 0  &  \ldots& 0 & -1 &-1\\
-1 & 3 & - 1 & 0 &    0&  0& \ddots& \vdots\\
0 &-1 & 3 & - 1 & 0 &    0&  \ddots&  \vdots\\
\vdots &0 &-1 & 3 & - 1 & 0 &  \ddots&  -1\\
\vdots & \ddots &\ddots & \ddots &  \ddots &  \ddots&  \ddots & -1\\
0 & \ddots & 0 & 0 &  -1 &  3&  -1 & -1\\
-1 & \ddots  &0 & 0 & 0 &  -1 &  3&  -1 \\[1.7mm]
-1  &\dots& \dots & -1 &-1& -1 &  -1 & k-1
\end{array}\right].\]

\begin{figure}
\[\begin{array}{c}
\begin{tikzpicture}[line cap=round,line join=round,>=triangle 45,x=1.0cm,y=1.0cm,scale = 1.5]
\draw [line width=.8pt] (5.,3.)-- (4.191784354657935,3.848009081626107);
\draw [line width=.8pt] (4.39006772869402,1.9998427342084801)-- (3.8288722643919675,3.028138208643161);
\draw [line width=.8pt] (4.39006772869402,1.9998427342084801)-- (5.5611954643020525,1.9717045255653174);

\draw [line width=.8pt] (4.39006772869402,1.9998427342084801)-- (5.,3.);
\draw [line width=.8pt] (5.,3.)-- (3.8288722643919675,3.028138208643161);
\draw [line width=.8pt] (3.8288722643919675,3.028138208643161)-- (4.191784354657935,3.848009081626107);
\draw [line width=.8pt] (4.191784354657935,3.848009081626107)-- (5.,3.);
\draw [line width=.8pt] (5.,3.)-- (5.028138208643161,4.1711277356080325);
\draw [line width=.8pt] (5.028138208643161,4.1711277356080325)-- (4.191784354657935,3.848009081626107);
\draw [line width=.8pt] (5.,3.)-- (5.848009081626107,3.8082156453420652);
\draw [line width=.8pt] (5.848009081626107,3.8082156453420652)-- (6.171127735608033,2.9718617913568393);
\draw [line width=.8pt] (6.171127735608033,2.9718617913568393)-- (5.,3.);
\draw [line width=.8pt] (5.,3.)-- (5.5611954643020525,1.9717045255653174);
\draw [line width=.8pt] (5.5611954643020525,1.9717045255653174)-- (6.171127735608033,2.9718617913568393);
\begin{scriptsize}
\draw [fill=black] (5.,3.) circle (1.5pt);
\draw[color=black] (5,2.7) node {$k$};
\draw [fill=black] (4.191784354657935,3.848009081626107) circle (1.5pt);
\draw[color=black] (4.05,4.05) node {$3$};
\draw [fill=black] (3.8288722643919675,3.028138208643161) circle (1.5pt);
\draw[color=black] (3.5,3.08) node {$2$};
\draw [fill=black] (5.028138208643161,4.1711277356080325) circle (1.5pt);
\draw[color=black] (5,4.381775026148471) node {$4$};
\draw [fill=black] (5.848009081626107,3.8082156453420652) circle (1.5pt);
\draw[color=black] (6.1,4) node {$k-3$};
\draw [fill=black] (6.171127735608033,2.9718617913568393) circle (1.5pt);
\draw[color=black] (6.62,3.05) node {$k-2$};
\draw [fill=black] (4.39006772869402,1.9998427342084801) circle (1.5pt);
\draw [fill=black] (5.5611954643020525,1.9717045255653174) circle (1.5pt);
\draw[color=black] (5.25,1.8) node {$k-1$};
\draw[color=black] (4.2,1.8) node {$1$};
\draw [fill=black] (5.438073645134634,3.989671690475049) circle (.9pt);
\draw [fill=black] (5.233105926888897,4.08039971304154) circle (.9pt);
\draw [fill=black] (5.64304136338037,3.898943667908557) circle (.9pt);
\end{scriptsize}
\end{tikzpicture}
 \end{array}\]
\caption{The wheel graph on $k$ vertices}\label{fig:wheels}
\end{figure}
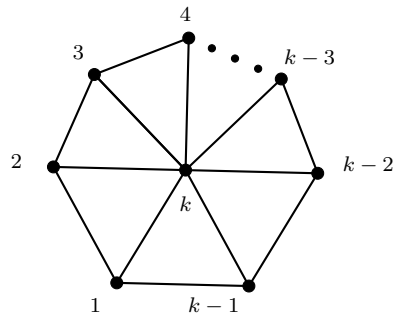
\begin{prop}\cite{BapatWheels}\label{pro:bapatwheels}
For $k \ge 2,$ the resistance distance between vertex $k$ and vertex $i$, $i\in\{1,\ldots, k-1\}$ in the wheel graph
is 
\[r(i,k)= \frac{F_{2k-2}^2}{F_{4k-4}-2F_{2k-2}}.\]
\end{prop}

Alternatively, we may formulate $r(i,k)$ in terms of recursive sequences as follows. $r(i,k)=\frac{N_k}{D_k}$ where $\{N_k\}_{k \ge 0}$
satisfies the recursion $G_n = 8 G_{n-1} - 8 G_{n-2} + G_{n-3} $ with initial conditions $\{G_0, G_1, G_2, G_3\} = 
\{1, 0, 1, 9\},$
and $\{D_k\}_{k \ge 0}$
satisfies the recursion
$G_n = 10 G_{n-1} - 23 G_{n-2}+10 G_{n-3} - G_{n-4}$ with   initial conditions
$\{G_0,\dotsc,G_4\} =
\{-1, 0, 1, 15, 128\}.$
\section{Conclusion}

This paper presents an introductory survey of resistance distance, a rich array of applications to which it is applied, a light review of computational methods, and some basic examples including  their associated recursions. It is hoped that this paper will encourage readers of the Quarterly to pursue this beautiful and active field.


\end{document}